\documentclass[11pt,letterpaper]{article}

\usepackage{sidecap}

\usepackage{graphics} 
\usepackage{epsfig} 
\usepackage{times} 
\usepackage{amsmath} 
\usepackage{amssymb}  
\usepackage{graphicx}
\usepackage{amssymb}
\usepackage[small]{caption}

\newcommand{\beqn}[1]{\begin{eqnarray}\label{#1}}
\newcommand{\eeqn}{\end{eqnarray}}
\newcommand{\beq}{\begin{eqnarray*}}
\newcommand{\eeq}{\end{eqnarray*}}
\newcommand{\bit}{\begin{itemize}}
\newcommand{\eit}{\end{itemize}}

\newcommand{\comment}[1]{}

\newcommand{\be}[1]{\begin{equation}\label{#1}}
\newcommand{\ee}{\end{equation}}

\newcommand{\bi}{\begin{itemize}}
\newcommand{\ei}{\end{itemize}}
\newcommand{\bc}{\begin{center}}
\newcommand{\ec}{\end{center}}

\newcommand{\ben}{\begin{enumerate}}
\newcommand{\een}{\end{enumerate}}

\title{Some Remarks on Input Choices for Biochemical Systems}

\author{Eduardo D. Sontag\\
Department of Mathematics and BioMaPS Institute\\
Rutgers University, New Brunswick, NJ 08903\\
\texttt{sontag@math.rutgers.edu}}

\begin{document}

\maketitle

Control systems theory concerns the study of open dynamical systems that
process time-dependent input signals (stimuli, ligands, controls, forcing
functions, test signals) into output signals (responses, measurements,
read-outs, reporters).  Such input/output (i/o) systems may be studied by
themselves, or as components (subsystems, modules) of larger systems.
Two extreme approaches to modeling systems (or particular subsystems) are as
follows~\cite{mct,icm94}:
\bi
\item
The \emph{black-box}, \emph{purely input/output}, or
\emph{phenomenological} description ignores all mechanistic detail and
characterizes behavior solely in terms of behavioral or stimulus-response
data.  This 
data may be summarized for example by transfer functions (for linear systems
and certain types of ``bilinear'' systems), Volterra, generating series, or
Wiener expansions (for general nonlinear systems), statistical models
that correlate inputs and outputs, neural network and other approximate models,
or machine learning prediction formalisms.
\item
The \emph{state-space}, \emph{internal}, or \emph{mechanistic} description,
in which all relevant players (``state variables'' such as proteins, mRNA, and
metabolites) are specified, and the forms for all reactions and reaction
constants are described.
\ei
In between these two extremes, one finds combinations (``gray box'' view).
Available mechanistic information is incorporated into the model, and this is
supplemented by i/o data, which serves to impose constraints on (or, in a
Bayesian approach, to provide prior distributions of) unknown parameters to
be further identified.
In engineering applications, black-box modeling is common, for
example, in chemical process control, and state-space modeling is
often used in mechanical and aerospace engineering.
Many biological models involve combinations of both approaches.

The subfield of \emph{realization theory} (``reverse-engineering'') deals with
the question of how much of the internal system can be deduced from its
input/output or black-box behavior.  One may pose questions about the
structure itself as well as about parameters, once that a concrete systems
structure has been hypothesized.
The basic results classify symmetry groups of possible internal realizations
of a given behavior, typically through linear algebra, algebraic-geometric,
and differential-geometric techniques, depending on the class of systems being
studied, or provide conditions for parameter ``identifiability'' from i/o
data.
(The term is not to be confused with \emph{identification theory}, the field
that incorporates stochastic aspects so as to deal with noisy data.  The
underlying theoretical questions of what is ultimately achievable, even in
ideal noise-free conditions, are easier to understand in a deterministic
setting.)

Key concepts used in realization theory are \emph{observability} (how much
information about internal dynamic variables, such as chemical concentrations,
can be ultimately obtained from time-varying stimulus/response
measurements?{}) and \emph{reachability (or controllability)} (what state
configurations can be attained, from a given initial state by manipulating
inputs, also called ``open-loop'' controls?{}).
Minimal models explaining observed input/output behavior typically
are characterized in terms of reachability and observability.
Also, identifiability of parameters can be reduced to an observability problem.
Many tests exist for observability and controllability, for linear and more
generally for large classes of nonlinear systems~\cite{mct}, although the
theory of nonlinear controllability is still very incomplete.

There are many other important subfields of control theory, among them the
study of \emph{disturbances and robustness} (how does one quantify the effect
of actuator and sensor disturbances, and the robustness to unstructured
perturbations?{}), \emph{feedback} (how does one design feedback laws or
``closed-loop controls'' so as to enhance desirable system characteristics
such as stability and robustness to structural perturbations, noise, and other
disturbances?{}), \emph{optimal control} (among those open-loops controls
which serve to drive one internal configuration of state variables to another,
which ones are optimal in appropriate senses?{}), and the question of hybrid
or \emph{discrete/continuous interfaces} (how do digital sensors, actuators,
limited-bandwidth communication channels, and feedback devices interact with
systems described by analog real-valued signals?{}).
Each of these subfields brings up its own set of challenges special to
systems biology~\cite{sysbio04EDS,ejc05EDS}.  

Here, we make some elementary remarks regarding how the choice of inputs being
applied to a system (pulses, steps, etc) affects parameter identifiability.
To keep the exposition accessible, we work out some extremely simple explicit
examples, instead of discussing general theoretical matters.  Some of the
material, such as that covering linear systems, is standard~\cite{mct}.
The nonlinear examples illustrate ideas from a paper in preparation.

\section*{Various Input Classes}

We discuss next how the choice of input signals impacts the \emph{theoretical}
possibility of identifying system responses (black-box behavior) and/or
internal structure (state-space view). 

Normally, the analysis of responses to small signals (or small perturbations
from steady states, leading to ``weakly activated'' systems) can be handled by
linear control theory, which provides a rich toolkit for the analysis of
biochemical networks.  On the other hand, if relevant signals are large, or if
nonlinear effects, such as receptor desensitization or dimerization, affect
signaling at a comparable time-scale, then linear tools do not suffice and new
techniques must be brought to bear.  Nonlinear control systems analysis has
undergone an exceptional period of development and maturation within the last
20 years, largely driven by application needs from aerospace, mechanical, and
other engineering disciplines.  The tools that have been developed reflect the
requirements from these fields, hence the emphasis on studying certain classes
of systems, such as Hamiltonian dynamics, which are less relevant to
biochemical networks.  Although many of the fundamental ideas are still
applicable to biology, much further theoretical work is needed.

\subsection*{Small-signal analysis: equivalence of step, impulse, and pulse responses}

When probing systems with signals of small magnitude, the full power of
the well-developed \emph{linear control theory} can be taken advantage of.
To illustrate the problems that arise in this context, we use a very simple
example, that of a single species $X$ whose concentration $x(t)$ is described
by an ordinary differential equation:
\[
\frac{dx}{dt}(t) = -ax(t) + bu(t)
\]
in which the variable $u(t)$ measures the, generally time varying,
concentration of an external input $U$, and $a,b$ are unknown positive
constants. 
(From now on, we will drop the time arguments $t$ when clear from the
context.)

For example, $X$ might represent the concentration of bound receptors,
and $U$ the concentration of a ligand which is present only in very small
quantities. 
In this case, $a$ is the degradation rate of $X$, due to internalization,
ligand unbinding, and other effects, and $b$ represents the concentration of
unoccupied receptors, which we consider as constant by assuming that their
abundance is large compared to the ligand.  
(If the concentration is not constant, then a linear model is not appropriate.)
Generally, variables such as $x$ and $u$ represent not absolute concentrations
but rather perturbations from a steady state.

We assume that the reporter being used provides a measurement
$y(t) = cx(t)$ equal to the product of the concentration of $X$ and
some (unknown) positive constant $c$.
Mechanistically, the form of $y$ might arise from an additional differential
equation:
\[
\varepsilon \frac{dy}{dt} = -y + cx
\]
which models the formation and decay of the reporter $Y$ (we took a unit decay
to simplify) and $0<\varepsilon \ll1$ indicates that this process happens at a much
faster time-scale than the production of $x$.  Using a quasi-steady state
approximation (singular perturbation), we set $\varepsilon =0$ and consider $y=cx$ as
the output.  Another way in which an output of this type could arise is as a
model of unknown measurement device calibration.

For any input $u(t)$ for $t>0$, and $x(0)=0$, we have, for $t>0$:
\[
y(t) = (k*u)(t) = \int_0^t k(t-r) u(r)\, dr
\]
(star denotes convolution), where
\[
k(r) = ce^{-at}b \,.
\]
The function $k(t)$ is called the \emph{impulse-response} of the system,
because a unit impulse $u=\delta _0$ at time $0$ (or, in practice, a very narrow
pulse with unit area) elicits the response $y(t)=k(t)$.
The function $K(t)$ is called the \emph{step-response}, because for any
step input $u(t)\equiv u_0$ for $t>0$, the output is $y(t) = K(t)u_0$ for $t>0$.

\emph{Knowledge of $K$ is equivalent to knowledge of $k$}, because
$k=K'$ and $K=\int_0^t k$.
Therefore, knowing the time-series response to \emph{any single step input} is
sufficient to determine $K(t)$ (just divide: $K(t) = y(t)/u_0$), which in turn
(convolution formula) \emph{determines the response to any other possible
input}, for instance pulses.  

In fact, knowing the response to basically \emph{any single input} is
theoretically enough for finding $k$.  A proof using Laplace transforms, is
clear from the fact that $\hat y(s) = W(s) \hat u(s)$, where $W(s)$ is the
frequency response of the system, that is to say, the Laplace transform of
$k$.  Thus, as long as $\hat u(s)$ does not vanish, one can solve for $W=\hat
y / \hat u$, and $k$ can then be obtained by Laplace inversion.
In practice, of course, one does not perform just a single input/output
experiment, since measurements will be noisy, nor does one obtain $k$ as $K'$,
since differentiation of an estimated step function is numerically difficult.
Several measurements are averaged, through a least-squares regression or
another technique.  This is discussed at length in textbooks in identification
theory.  It is easier, however, to explain the theory in the noise-free case.

\subsection*{Steady state measurements do not provide enough information about
    transient behavior}

Steady-state experiments measure the steady state output $y_{\rm ss}$ of the
reporter variable as a function of the magnitude of a step input.
With the above example, $y_{\rm ss} = c x_{\rm ss} = \frac{cb}{a}u_0$.
In other words, the \emph{only} information about the three parameters $a,b,c$
that can be obtained from such experiments is the combination $\gamma =bc/a$
(called the steady state gain of the system).
For example, the two systems (i) $dx/dt=-x+u$ with $y=x$ and 
(ii) $dx/dt=-2x+u$ with $y=2x$ have the same steady-state gain, $\gamma =1$.
The i/o behaviors of these two systems are different, however: when using a
step input $u(t)\equiv 1$, and initial condition $x(0)=0$,
the output of the first system is $y(t)=1-e^{-t}$, but
the output of the second system is a different function, $y(t)=1-e^{-2t}$.
\emph{Only transient measurements reveal this distinction.}

\subsection*{Parameter unidentifiability leads to equivalence classes of
    parameters} 

Still with the same example, we note that the product $bc$ is an i/o invariant
of the system, because $b$ and $c$ appear together as a product in the
convolution formula.
No matter what input is applied (assuming $x(0)=0$), two systems for which the
same product $bc$ is the same will give the same output.
No amount of i/o experimentation will permit finding $b$ and $c$ individually.

On the other hand, the decay $a$ can be determined uniquely from i/o
experiments.  For example, with a unit step input, we may measure
$y(t)=K(t)$, from which we may evaluate $k(t)$; hence $k(0)=K'(0)=cb$ and
$k'(0)=K''(0)=-cab$ are known, and $a = -K''(0)/K'(0)$ can be estimated.
One says that the two parameter triples $(a_1,b_1,c_1)$ and $(a_2,b_2,c_2)$
are \emph{equivalent} if they cannot be distinguished by any i/o experiment.
Thus, we showed that two triples are equivalent if and only if $b_1c_1=b_2c_2$
and $a_1=a_2$.
Another way to say state this fact is through a one-parameter group of
symmetries:
$(a,b,c)$ is equivalent to $(a,Tb,cT^{-1})$ for all $T\not= 0$.
(These symmetries amount to rescaling the concentration $x(t)$: in different
units, one has a new variable $z(t)=Tx(t)$, which satisfies the differential
equation $dz/dt=-az+Tbu$, and the reporter is now and $y=cT^{-1}x$.)

A far-reaching generalization of this result is available for linear systems
of arbitrary dimension $n$ (number of species), not merely $n=1$ as in our
example: systems $dx/dt=Ax+bu$, $y=cx$, where now $A$ is an $n\times n$ matrix, $b$
is a column $n$-vector, and $c$ is row $n$-vector.  Under genericity conditions
(``reachability and observability''), the following theorem is true~\cite{mct}:
two triples $(A_i,b_i,c_i)$, $i=1,2$ are i/o equivalent if and only if
$(TA_1T^{-1},Tb_1,cT_1^{-1})=(A_2,b_2,c_2)$ for some invertible matrix $T$.
(The theorem is also true for multivariable inputs and outputs, in which case
both $B$ and $C$ are matrices.)

Continuing with our example, we may think of our system as a ``gray box'' in
the sense that we do not know the complete mechanism (all the parameters) but
neither do we have complete ignorance, since we are able to measure $a$ and
the product $bc$.
Now suppose that a new experimental design allows us not merely to measure
$y(t)$ but also do directly measure the parameter $c$ (for example, the
binding of $X$ to species $Y$ may be characterized by means of a biochemical
experiment).  Now it is possible to also determine $b$, since $b=K'(0)/c$.
Or perhaps, instead of a biochemical experiment, a new probe is devised, which
allows one to directly measure the concentration $x(t)$.  Then, knowing that
$dx/dt=-ax+bu$, and from knowledge of $x(t)$ and $x'(t)$ at two values $t=t_i$,
one can determine both $a$ and $b$, and hence, $c=K'(0)/b$ as well.
In either case, we have gone from a ``gray box'' description to a complete
mechanistic description.
Partial knowledge of $b$ or $c$ can also be understood in these terms: suppose
that we know that $K'(0)=1$ and that an imprecise biochemical assay provides
the range $c\in [0.01,0.1]$.  Then, we know that $b\in [10,100]$, a ``less gray''
level of description than total ignorance of $b$.

These ideas are particularly useful in a Bayesian context.  The i/o knowledge
given by the value of $a$ and the product $bc$, combined with physically
realistic ranges for the unknown $b$ and $c$, can be used to provide priors 
for parameter estimation through further experiments (or to evaluate
the reliability of using parameters from a homologous protein, for example).
That is, the invariants provide a ``prior belief'' about the parameters
$\pi =(a,b,c)$, summarized as a probability $P(\pi )$.  After an experiment
$e$, this probability may then be updated to obtain a posterior distribution
$P(\pi |e)$ through an application of Bayes Rule: 
\[
P(\pi |e) = \frac{P(e|\pi )P(\pi )}{P(e)}, 
\]
where $P(e|\pi )$ is the likelihood of the experimental results under the
assumed parameters $\pi $.


\subsection*{Step inputs are not sufficient for nonlinear identification}

There are strong theoretical justifications for the need for rich input
classes in the context of reverse engineering biological systems.
We use a simplified, and hence somewhat unrealistic, example in order to
explicitly illustrate this point, but the phenomenon scales to larger systems.

Consider the following system, which describes the time-evolution of the
concentrations $x(t)$ and $z(t)$ of two chemical species $X$ and $Z$:
\beq
\frac{dx}{dt} &=& -\lambda x + u^2\\
\frac{dz}{dt} &=& -\lambda z + u
\eeq
with initial concentrations $x(0)=0$ and $z(0)=0$, 
and where $u(t)$ is the (possibly time-varying) concentration of an input (for
instance a ligand) $U$ that helps promote formation of $X$ and $Z$.
Our goal is to estimate the unknown parameter $\lambda $, the degradation rate of
$X$ and $Z$.
Note that as long as no input $U$ is applied ($u(t)=0$), the concentrations of
$X$ and $Z$ remain at steady state (zero).

We assume that the only available measurement is the (time-varying)
concentration $y(t)$ of a reporter $Y$, which we take as a function of the
instantaneous concentrations of $X$ and $Z$ according to the formula:
\[
y = x + u(a-z)
\]
were $a$ is some known positive constant.
(One could view $y$ as obtained through a quasi-steady state approximation
applied to 
\[
\varepsilon \frac{dy}{dt} = -y + x + u(a-z)
\]
with $0<\varepsilon \ll1$, which describes a reaction in which a certain substance
$A$ dimerizes with $U$ in order to promote formation of $y$, but $A$ is
sequestered by $Z$ with 1-1 stoichiometry, so the free amount of $A$
equals $a-z=A_{\rm total}-A_{\rm bound}$, because the concentration of 
$A_{\rm bound}$ equals the concentration of $Z$.)

%
%
%
%
%

If we use a step function input, $u(t)\equiv u_0$ for $t>0$, then
$x(t)\equiv z(t)u_0$, and therefore $y(t)\equiv au_0$, independently of the
actual value of $\lambda $.
Thus, \emph{with step inputs}, no matter how many experiments (different
values of $u_0$) are carried out, \emph{it is impossible to obtain any
information about $\lambda $}.

On the other hand, suppose that we use a unit ramp: $u(t)=t$.
Then $y''''(t)= -\lambda  (\lambda  t - 2) e^{-\lambda  t}$, and so
\[
\lim_{t\searrow 0}y''''(t)=2\lambda 
\]
and therefore $\lambda $ can be identified by differentiation.
(Of course, differentiation will not be used in practice, nor is it needed
since one may use other techniques for estimating $\lambda $, but we are just
proving theoretical identifiability.)
In summary, a single ramp input is enough to identify $\lambda $, but no amount of
even complete transient time-measurements of step responses will achieve
this result.

Instead of ramps, we could also use pulses.  Suppose that we consider a unit
pulse: $u(t)=1$ for $0\leq t\leq 1$ and $u(t)=0$ for $t>1$.
Now, $x(t)=z(t)=\frac{1}{\lambda }(1-e^{-\lambda t})$ for $t\in [0,1]$ and
$x(t)=z(t)=\frac{1}{\lambda }(e^{\lambda }-1)e^{-\lambda t}$ for $t\geq 1$, so
$y(t)\equiv x(t)+0(A-z(t)) = x(t)$ for $t\geq 1$, which means that
\[
\lim_{t\searrow 1}y'(t)=e^{-\lambda }-1 \,,
\]
and therefore $\lambda $ can be calculated as:
\[
\lambda =-\ln\left(1+\lim_{t\searrow 1}y'(t)\right)\,.
\]
Once again, a non-step signal is successful in identifying the parameter $\lambda $,
but steps are not enough.

Similar examples can be given with more parameters and more complicated
systems.  Also, we took the same degradation rates for $X$ and $Z$, because we
are considering a noise-free ideal situation.  If these degradation rates
would not necessarily be equal but are close, then $y(t)$ will depend on $\lambda $,
even in the step case.  However, the dependence will be through a very small
multiplier coefficient, indicating a very small sensitivity and therefore very
high estimator variance (Cramer-Rao bounds), thus rendering estimation
practically impossible.  It is possible to prove (paper in preparation)
that for bilinear systems (those in which the only possible products are
between inputs and states), pulses provide as much theoretical information as
any possible time-varying input signals. 

\section*{Remarks}

We have only discussed \emph{parameter} identification, as opposed to
structure identification.  For small-signal (linear) analysis, the two
problems are the same (at least if an upper bound on the number of species in
the system is known), since zero entries in the differential equations
indicate non-existent interactions.  For nonlinear systems, the problem of
identification of structure is far more difficult.  It is the subject (at
least for the ideal noise-free case) of nonlinear realization
theory~\cite{icm94}.
The parameter estimation problem is a special case of a nonlinear
observability problem (thinking of unknown parameters as constant states), and
a rich theory allows the use of Lie-algebra and differential-algebra
techniques for determining observability~\cite{mct}.

There are results available characterizing the number of experiments
required for identifying a given number of parameters in a system~\cite{ident_experiments02}.
Another relevant theoretical mathematical fact is that generic ``random''
single inputs are enough for identification~\cite{icm94}, and this fact
is very useful in the context of identification of certain biochemical
systems~\cite{vanriel06}.


\begin{thebibliography}{1}

\bibitem{mct}
E.D. Sontag.
\newblock {\em Mathematical control theory: deterministic finite dimensional
  systems (2nd ed.)}.
\newblock Springer-Verlag New York, Inc., New York, NY, USA, 1998.

\bibitem{icm94}
E.D. Sontag.
\newblock Spaces of observables in nonlinear control.
\newblock In {\em Proceedings of the International Congress of Mathematicians,
  Vol.\ 1, 2 (Z\"urich, 1994)}, pages 1532--1545, Basel, 1995. Birkh\"auser.

\bibitem{sysbio04EDS}
E.D. Sontag.
\newblock Some new directions in control theory inspired by systems biology.
\newblock {\em IEE Proc. Systems Biology}, 1:9--18, 2004.

\bibitem{ejc05EDS}
E.D. Sontag.
\newblock Molecular systems biology and control.
\newblock {\em Eur. J. Control}, 11(4-5):396--435, 2005.

\bibitem{ident_experiments02}
E.D. Sontag.
\newblock For differential equations with {$r$} parameters, {$2r+1$}
  experiments are enough for identification.
\newblock {\em J. Nonlinear Sci.}, 12(6):553--583, 2002.

\bibitem{vanriel06}
N.A.W. van Riel and E.D. Sontag.
\newblock Parameter estimation in models combining signal transduction and
  metabolic pathways: The dependent input approach.
\newblock {\em IEE Proc.\ Systems Biology}, to appear, 2006.

\end{thebibliography}

\end{document}